\documentclass[a4paper,11pt]{amsart}

\usepackage[dvips]{graphicx}
\usepackage{amsmath,amssymb,amsfonts,amsthm}
\newtheorem{theorem}{Theorem}[section]
\newtheorem{maintheorem}{Theorem}

\newtheorem{secondtheorem}{Theorem}

\newtheorem{lemma}[theorem]{Lemma}
\theoremstyle{definition}
\newtheorem{definition}[theorem]{Definition}
\newtheorem{example}[theorem]{Example}
\newtheorem{remark}[theorem]{Remark}
\numberwithin{equation}{section}

\begin{document}
\title[Normal Forms of Levi-Flat hypersurfaces]{On Normal Forms of Singular Levi-Flat Real Analytic Hypersurfaces}
\author{Arturo Fern\' andez-P\' erez}
\address{Instituto de Matem\' atica Pura e Aplicada, IMPA}
\curraddr{Estrada Dona Castorina, 110, 22460-320. Rio de Janeiro, RJ, Brazil.}
\email{afernan@impa.br}

\subjclass[2010]{Primary 32V40 - 37F75}
\keywords{Levi-flat hypersurfaces - Holomorphic foliations}

\begin{abstract}
Let $F(z)=\mathcal{R}e(P(z)) + h.o.t$ be such that $M=(F=0)$ defines a germ of real analytic Levi-flat at $0\in\mathbb{C}^{n}$, $n\geq{2}$, where $P(z)$ is a homogeneous polynomial of degree $k$ with an isolated singularity at $0\in\mathbb{C}^{n}$ and Milnor number $\mu$. We prove that there exists a holomorphic change of coordinate $\phi$ such that $\phi(M)=(\mathcal{R}e(h)=0)$, where $h(z)$ is a polynomial of degree $\mu+1$ and $j^{k}_{0}(h)=P$.
\end{abstract}
\maketitle
\section{Introduction and Statement of the results}
\par Let $M$ be a germ at  $0\in\mathbb{C}^{n}$ of a real codimension one irreducible analytic set. For the sake of simplicity we will denote germs and representative of germs by the same letter. Since $M$ is real analytic of codimension one and irreducible, it can be defined in $\mathbb{C}^{n}$ by $(F=0)$, where $F$ is an irreducible germ of real analytic function.
The singular set of $M$ is defined  by $sing(M)=(F=0)\cap(dF=0)$  and its smooth part $(F=0)\backslash(dF=0)$ will denoted by $M^{*}$. The Levi distribution $L$ on $M^{*}$ is defined by $L_{p}:=ker(\partial{F}(p))\subset T_{p}M^{*}=ker(dF(p))$, for any $p\in M^{*}$.
\begin{definition}
We say $M$ is Levi-flat if the Levi distribution on $M^{*}$ is integrable.
\end{definition}
\begin{remark} The integrability condition of $L$ implies that $M^{*}$ is tangent to a real codimension one foliation $\mathcal{L}$. Since the hyperplanes $L_{p}$, $p\in M^{*}$, are complex, the leaves of $\mathcal{L}$ are complex codimension one holomorphic submanifolds immersed on $M^{*}$.
\end{remark}
\begin{remark}
If the hypersurface $M$ is defined by $(F=0)$ then the Levi distribution $L$ on $M$ can be defined by the real analytic 1-form $\eta=i(\partial{F}-\bar{\partial}F)$, which will be called the Levi 1-form of $F$. The integrability condition is equivalent to $(\partial{F}-\bar{\partial}F)\wedge\partial\bar{\partial}F|_{M^{*}}=0$
\end{remark}

\par In the case of a real analytic smooth Levi-flat hypersurface $M$ in $\mathbb{C}^{n}$, its local structure is very well understood, according to E. Cartan, around each $p\in M$ we can find local holomorphic coordinates $z_{1},\ldots,z_{n}$ such that $M=\{\mathcal{R}e(z_{1})=0\}$. 
\par More recently D. Burns and X. Gong [B-G] have proved an analogous result in the case $M=F^{-1}(0)$ Levi-flat, where $F:(\mathbb{C}^{n},0)\rightarrow(\mathbb{R},0)$, $n\geq 2$, is a germ of real analytic function such that  $$F(z_{1},\ldots,z_{n})=\mathcal{R}e(z_{1}^{2}+\ldots+z_{n}^{2})+h.o.t.$$ 
\par They show that there exists a germ of biholomorphism $\phi:(\mathbb{C}^{n},0)\rightarrow(\mathbb{C}^{n},0)$ such that $\phi(M)=(\mathcal{R}e(z_{1}^{2}+\ldots+z_{n}^{2})=0)$. 
\par In [C-L], the authors prove the above result by using the theory of holomorphic foliations. In this paper we are interested in finding similar normal forms in a situation more general. Our main result is the following:
\begin{maintheorem}\label{main-theorem}
Let $M=F^{-1}(0)$, where $F:(\mathbb{C}^{n},0)\rightarrow(\mathbb{R},0)$, $n\geq 2$, be a germ of irreducible real analytic function such that
\begin{enumerate}
\item $F(z_{1},\ldots,z_{n})=\mathcal{R}e(P(z_{1},\ldots,z_{n}))+h.o.t$,
where $P$ is a homogeneous polynomial of degree $k$ with an isolated singularity at $0\in\mathbb{C}^{n}$.
\item The Milnor number of $P$ at $0\in\mathbb{C}^{n}$ is $\mu$.
\item $M$ is Levi-flat.
\end{enumerate}
Then there exists a germ of biholomorphism $\phi:(\mathbb{C}^{n},0)\rightarrow(\mathbb{C}^{n},0)$ such that $\phi(M)=(\mathcal{R}e(h)=0)$, where $h(z)$ is a polynomial of degree $\mu+1$ and $j^{k}_{0}(h)=P$. 
\end{maintheorem}

\begin{remark}
In particular, we obtain the result of [B-G].
\end{remark}
\begin{secondtheorem}\label{second-theorem}
In the same spirit we have the following generalization: Let $M=F^{-1}(0)$, where $F:(\mathbb{C}^{n},0)\rightarrow(\mathbb{R},0)$, $n\geq 3$, be a germ of irreducible real analytic function such that
\begin{enumerate}
\item $F(z_{1},\ldots,z_{n})=\mathcal{R}e(Q(z_{1},\ldots,z_{n}))+h.o.t$,
where $Q$ is a quasihomogeneous polynomial of degree $d$ with an isolated singularity at $0\in\mathbb{C}^{n}$.
\item $M$ is Levi-flat.
\end{enumerate}
Then there exists a germ of biholomorphism $\phi:(\mathbb{C}^{n},0)\rightarrow(\mathbb{C}^{n},0)$ such that $$\phi(M)=(\mathcal{R}e(Q(z)+\sum_{j} c_{j}e_{j}(z))=0),$$ where $e_{1},\ldots,e_{s}$ are the elements of the monomial basis of the local algebra of $Q$ such that $deg(e_{j})>d$
and $c_{j}\in\mathbb{C}$.
\end{secondtheorem}

\noindent {\bf Acknowledgments}

\par I want to express my thanks to my advisor
 Alcides Lins Neto, for valuable conversations. The author would like to thank also IMPA, where this work was developed. I also want to thank the referee for his suggestions.

\section{Preliminaries}
Let us fix some notations that will be used from now on.
\begin{enumerate}
\item $\mathcal{O}_{n}$ : The ring of germs of holomorphic functions at $0\in\mathbb{C}^{n}$.\\
 $\mathcal{O}(U)=$ set of holomorphic functions in the open set $U\subset\mathbb{C}^{n}$. 
\item $\mathcal{O}^{*}_{n}=\{f\in\mathcal{O}_{n}/f(0)\neq 0\}$. \\
$\mathcal{O}^{*}(U)=\{f\in\mathcal{O}(U)/f(z)\neq 0,\forall z\in U\}$.  
\item $\mathcal{M}_{n}=\{f\in\mathcal{O}_{n}/f(0)=0\}$ maximal ideal of $\mathcal{O}_{n}$.
\item $\mathcal{A}_{n}$ : The ring of germs at $0\in\mathbb{C}^{n}$ of complex 
valued real analytic functions.

\item $\mathcal{A}_{n \mathbb{R}}$ : The ring of germs at $0\in\mathbb{C}^{n}$ of real valued real analytic functions. Note that $F\in\mathcal{A}_{n}$ is in 
$\mathcal{A}_{n \mathbb{R}}$ if and only if $F=\bar{F}$.

\item $Diff(\mathbb{C}^{n},0)$ : The group of germs at $0\in\mathbb{C}^{n}$ of holomorphic diffeomorphisms $f:(\mathbb{C}^{n},0)\rightarrow(\mathbb{C}^{n},0)$ with the operation of composition.
\item $j^{k}_{0}(f)$ : The $k$-jet at $0\in\mathbb{C}^{n}$ of $f\in\mathcal{O}_{n}$.

\end{enumerate}

\begin{definition}
 Two germs $f,g\in\mathcal{O}_{n}$ are said to be right equivalent, if there exists $\phi\in Diff(\mathbb{C}^{n},0)$ such that $f\circ\phi^{-1}=g$.
\end{definition}

The local algebra of $f\in\mathcal{O}_{n}$ is by definition $$A_{f}=\mathcal{O}_{n}/(\partial{f}/\partial{z}_{1},\ldots,\partial{f}/\partial{z}_{n}).$$
\begin{definition}
Define by $\mu(f,0):=dim A_{f}$, the Milnor number of $f$ at $0\in\mathbb{C}^{n}$.
\end{definition}

\par Morse Lemma can now be rephrased by saying that if $0\in\mathbb{C}^{n}$ is an isolated singularity of $f$ with Milnor number $\mu(f,0)=1$ then $f$ is right equivalent to its second jet. The next lemma is a generalization of Morse's Lemma. We refer to [A-G-V], pg.121.
\begin{lemma}\label{morse-lemma}
Suppose $0\in\mathbb{C}^{n}$ is an isolated singularity of $f\in\mathcal{M}_{n}$ with Milnor number $\mu$. Then $f$ is right equivalent to $j^{\mu+1}_{0}(f)$.

\end{lemma} 

\subsection{The complexification} 
In this section we state some general facts about complexification of germs of real analytic functions. 
\par Given $F\in\mathcal{A}_{n}$; we can write its Taylor series at $0\in\mathbb{C}^{n}$ as
\begin{equation}\label{equa-uno}
F(z)=\sum_{\mu,\nu}F_{\mu\nu}z^{\mu}\bar{z}^{\nu},
\end{equation}
 where 
$F_{\mu \nu}\in\mathbb{C}$, $\mu=(\mu_{1},\ldots,\mu_{n})$, 
$\nu=(\nu_{1},\ldots,\nu_{n})$, $z^{\mu}=z_{1}^{\mu_{1}}\ldots z_{n}^{\mu_{n}}$,
$\bar{z}^{\nu}=\bar{z}_{1}^{\nu_{1}}\ldots \bar{z}_{n}^{\nu_{n}}$. When $F\in\mathcal{A}_{n\mathbb{R}}$ then the coefficients $F_{\mu \nu}$ satisfy
$$\bar{F}_{\mu \nu}=F_{\nu \mu}.$$
\par The complexification $F_{\mathbb{C}}\in\mathcal{O}_{2n}$ of $F$ is defined by the series 
\begin{equation}\label{equa-dos}
F_{\mathbb{C}}(z,w)=\sum_{\mu,\nu}F_{\mu\nu}z^{\mu}w^{\nu}.
\end{equation}
\par If $F\in\mathcal{A}_{n\mathbb{R}}$, $F(0)=0$ and $M=F^{-1}(0)$ defines a Levi-flat, the complexification $\eta_{\mathbb{C}}$ of its Levi 1-form $\eta=i(\partial{F}-\bar{\partial}F)$ can be written as 
$$\eta_{\mathbb{C}}=i(\partial_{z}F_{\mathbb{C}}-\partial_{w}F_{\mathbb{C}})=i
\sum_{\mu,\nu}(F_{\mu\nu}w^{\nu}d(z^{\mu})-F_{\mu\nu}z^{\mu}d(w^{\nu})).$$
\par The complexification $M_{\mathbb{C}}$ of $M$ is defined as $M_{\mathbb{C}}=F^{-1}_{\mathbb{C}}(0)$ and its smooth part is $M^{*}_{\mathbb{C}}=M_{\mathbb{C}}\backslash(dF_{\mathbb{C}}=0)$. The integrability condition of $\eta=i(\partial{F}-\bar{\partial}F)|_{M^{*}}$ implies that $\eta_{\mathbb{C}}|_{M^{*}_{\mathbb{C}}}$ is integrable. Therefore 
$\eta_{\mathbb{C}}|_{M^{*}_{\mathbb{C}}}=0$ defines a foliation $\mathcal{L}_{\mathbb{C}}$ on $M^{*}_{\mathbb{C}}$ that will be called the complexification of $\mathcal{L}$.
\begin{definition}
The algebraic dimension of $sing(M)$ is the complex dimension of the singular set of $M_{\mathbb{C}}$.
\end{definition}

\par Consider a germ at $0\in\mathbb{C}^{2}$ of real analytic Levi-flat $M=(F=0)$,  where $F$ is irreducible in $\mathcal{A}_{2\mathbb{R}}$. Let $F_{\mathbb{C}}$, $M_{\mathbb{C}}=(F_{\mathbb{C}}=0)\subset(\mathbb{C}^{4},0)$ and $M_{\mathbb{C}}^{*}$ be as before. We will assume that the power series that defines $F_{\mathbb{C}}$ converges in a neighborhood of $\bar{\bigtriangleup}=\{(z,w)\in\mathbb{C}^{4}/|z|, |w|\leq 1\}$, so that $F_{\mathbb{C}}(z,\bar{z})=F(z)$ for all $|z|\leq 1$.
\par Let $V:=M_{\mathbb{C}}^{*}\backslash  sing(\eta_{\mathbb{C}}|_{M^{*}_{\mathbb{C}}})$ and denote $L_{p}$ the leaf of $\mathcal{L}_{\mathbb{C}}$ through $p$, where $p\in V$. In this situation we have the following important Lemma of [C-L].
\begin{lemma}\label{closed-lema}

In the above situation, for any $p=(z_{0},w_{0})\in V$ the leaf $L_{p}$ is closed in $M^{*}_{\mathbb{C}}$.

\end{lemma}

\par In the proof of theorem 1 we will use the following result of [C-L].
\begin{theorem}\label{alcides-theorem}
 Let $M=F^{-1}(0)$ be a germ of an irreducible real analytic Levi-flat hypersurface at $0\in\mathbb{C}^{n}$, $n\geq{2}$, with Levi 1-form $\eta=i(\partial{F}-\bar{\partial} F)$. Assume that the algebraic dimension of $sing(M)\leq 2n-4$. Then there exists an unique germ at $0\in\mathbb{C}^{n}$ of holomorphic codimension one foliation $\mathcal{F}_{M}$ tangent to $M$, if one of the following conditions is fulfilled:
\begin{enumerate}
\item $n\geq 3$ and $cod_{M_{\mathbb{C}}^{*}}(sing(\eta_{\mathbb{C}}|_{M_{\mathbb{C}}^{*}}))\geq 3$.
\item  $n\geq 2$, $cod_{M_{\mathbb{C}}^{*}}(sing(\eta_{\mathbb{C}}|_{M_{\mathbb{C}}^{*}}))\geq 2$ 
and $\mathcal{L}_{\mathbb{C}}$ has a non-constant holomorphic first integral.
\end{enumerate}
Moreover, in both cases the foliation $\mathcal{F}_{M}$ has a non-constant holomorphic first integral $f$ such that $M=(Re(f)=0)$.
\end{theorem}

\section{Proof of theorem 1} 

\par Let $M=F^{-1}(0)\subset(\mathbb{C}^{n},0)$ be a Levi-flat, where $F(z)=\mathcal{R}e(P(z))+h.o.t$ with $P$ be a homogeneous polynomial of degree $k\geq 2$ with an isolated singularity at $0\in\mathbb{C}^{n}$ and Milnor number $\mu$. We want to prove that there exists $\phi\in Diff(\mathbb{C}^{n},0)$ such that $\phi(M)=(\mathcal{R}e(h)=0)$, where $h$ is a 
polynomial of degree $\mu+1$.
\par The idea is to use theorem \ref{alcides-theorem} to prove that there exists a germ $f\in\mathcal{O}_{n}$ such that the foliation $\mathcal{F}$ defined by $df=0$ is tangent to $M$ and $M=(\mathcal{R}e(f)=0)$. The foliation $\mathcal{F}$ can viewed as an extension to a neighborhood of $0\in\mathbb{C}^{n}$ of the Levi foliation $\mathcal{L}$ on $M^{*}$.
\par Suppose for a moment that $M=(\mathcal{R}e(f)=0)$ and let us conclude the proof. Without lost of generality, we can suppose that $f$ is not a power in $\mathcal{O}_{n}$. In this case $\mathcal{R}e(f)$ is irreducible (cf.
[C-L]). This implies that $\mathcal{R}e(f)=U.F$, where $U\in\mathcal{A}_{n \mathbb{R}}$ and $U(0)\neq 0$. Let $\sum_{j\geq k}f_{j}$ be the taylor series of $f$, where $f_{j}$ is a homogeneous polynomial of degree $j$, $j\geq k$.
 Then 
$$\mathcal{R}e(f_{k})=j^{k}_{0}(\mathcal{R}e(f))=j^{k}_{0}(U.F)=U(0).\mathcal{R}e(P(z_{1},\ldots,z_{n})).$$
Hence $f_{k}(z_{1},\ldots,z_{n})=U(0).P(z_{1},\ldots,z_{n})$. We can suppose that $U(0)=1$, so that 
\begin{eqnarray}
f(z)=P(z)+h.o.t
\end{eqnarray}
In particular, $\mu=\mu(f,0)=\mu(P,0)$, $f\in\mathcal{M}_{n}$, because $P$ has isolated singularity at $0\in\mathbb{C}^{n}$. Hence by lemma \ref{morse-lemma}, $f$ is right equivalent to $j^{\mu+1}_{0}(f)$, i.e. there exists $\phi\in Diff(\mathbb{C}^{n},0)$ such that $h:=f\circ\phi^{-1}=j^{\mu+1}_{0}(f)$. Therefore, $\phi(M)=(\mathcal{R}e(h)=0)$ and this will conclude the proof of theorem 1.
\par Let us prove that we can apply theorem \ref{alcides-theorem}. We can write 
\begin{equation}
F(z)=\mathcal{R}e(P(z_{1},\ldots,z_{n}))+H(z_{1},\ldots,z_{n}),
\nonumber
 \end{equation}
where $H:(\mathbb{C}^{n},0)\rightarrow(\mathbb{R},0)$ is a germ of real-analytic function and $j^{k}_{0}(H)=0$. For simplicity, we assume that $P$ has real coefficients. Then we get the complexification

$$F_{\mathbb{C}}(z,w)=\frac{1}{2}(P(z)+P(w))+H_{\mathbb{C}}(z,w)$$
and $M_{\mathbb{C}}=F^{-1}_{\mathbb{C}}(0)\subset(\mathbb{C}^{2n},0)$. In the general case, replacing $P(w)=\sum a_{j}w^{j}$ by $\tilde{P}(w)=\sum\bar{a}_{j}w^{j}$, we will recover each step of proof.

\par Since $P(z)$ has an isolated singularity at $0\in\mathbb{C}^{n}$, we get $sing(M_{\mathbb{C}})=\{0\}$, and so the algebraic dimension of $sing(M)$ is $0$. On other hand, the complexification of $\eta=i(\partial{F}-\bar{\partial}F)$ is $$\eta_{\mathbb{C}}=i(\partial_{z}F_{\mathbb{C}}-\partial_{w}F_{\mathbb{C}}).$$

Recall that $\eta|_{M^{*}}$ and $\eta_{\mathbb{C}}|_{M^{*}_{\mathbb{C}}}$ define $\mathcal{L}$ and $\mathcal{L}_{\mathbb{C}}$. Now we compute $sing(\eta_{\mathbb{C}}|_{M^{*}_{\mathbb{C}}})$. We can write $dF_{\mathbb{C}}=\alpha+\beta$, with $$\alpha=\sum_{j=1}^{n}\frac{\partial{F_{\mathbb{C}}}}{\partial{z}_{j}}dz_{j}:=
\frac{1}{2}\sum_{j=1}^{n}(\frac{\partial{P}}{\partial{z}_{j}}(z)+A_{j})dz_{j}$$
and $$\beta=\sum_{j=1}^{n}\frac{\partial{F_{\mathbb{C}}}}{\partial{w}_{j}}dw_{j}:=
\frac{1}{2}\sum_{j=1}^{n}(\frac{\partial{P}}{\partial{w}_{j}}(w)+B_{j})dw_{j},$$
where $\frac{1}{2}\sum_{j=1}^{n}A_{j}dz_{j}=\sum_{j=1}^{n}\frac{\partial{H_{\mathbb{C}}}}{\partial{z}_{j}}dz_{j}$ and 
$\frac{1}{2}\sum_{j=1}^{n}B_{j}dw_{j}=\sum_{j=1}^{n}\frac{\partial{H_{\mathbb{C}}}}{\partial{w}_{j}}dw_{j}$.
\par Then $\eta_{\mathbb{C}}=i(\alpha-\beta)$, and so 
\begin{eqnarray}
\eta_{\mathbb{C}}|_{M^{*}_{\mathbb{C}}}=(\eta_{\mathbb{C}}+idF_{\mathbb{C}})|_{M^{*}_{\mathbb{C}}}=
2i\alpha|_{M^{*}_{\mathbb{C}}}=-2i\beta|_{M^{*}_{\mathbb{C}}}.
\end{eqnarray}

In particular, $\alpha|_{M^{*}_{\mathbb{C}}}$ and $\beta|_{M^{*}_{\mathbb{C}}}$
define $\mathcal{L}_{\mathbb{C}}$. Therefore  $sing(\eta_{\mathbb{C}}|_{M^{*}_{\mathbb{C}}})$ can be splited in two parts. Let $M_{1}=\{(z,w)\in M_{\mathbb{C}}| \frac{\partial{F}_{\mathbb{C}}}{\partial{w_{j}}}\neq 0$ for some $j=1,\ldots,n\}$ and $M_{2}=\{(z,w)\in M_{\mathbb{C}}| \frac{\partial{F}_{\mathbb{C}}}{\partial{z_{j}}}\neq 0$ for some $j=1,\ldots,n\}$, note that $M_{\mathbb{C}}=M_{1}\cup M_{2}$; if we denote by  
$$X_{1}:=M_{1}\cap\{\frac{\partial{P}}{\partial{z}_{1}}(z)+A_{1}=
\ldots=\frac{\partial{P}}{\partial{z}_{n}}(z)+A_{n}=0\}$$ and $$X_{2}:=M_{2}\cap\{\frac{\partial{P}}{\partial{w}_{1}}(w)+B_{1}=
\ldots=\frac{\partial{P}}{\partial{w}_{n}}(w)+B_{n}=0\},$$ 

then $sing(\eta_{\mathbb{C}}|_{M^{*}_{\mathbb{C}}})=X_{1}\cup X_{2}$. Since $P\in\mathbb{C}[z_{1},\ldots,z_{n}]$ has an isolated singularity at $0\in\mathbb{C}^{n}$, we  conclude that $cod_{M^{*}_{\mathbb{C}}}sing(\eta_{\mathbb{C}}|_{M^{*}_{\mathbb{C}}})=n$.

\par If $n\geq 3$,  we can directly apply Theorem \ref{alcides-theorem} and the proof ends.
In the case $n=2$, we are going to prove that $\mathcal{L}_{\mathbb{C}}$ has a non-constant holomorphic first integral.
\par We begin by a blow-up at $0\in\mathbb{C}^{4}$. Let $F(x,y)=\mathcal{R}e(P(x,y))+h.o.t$ and $M=F^{-1}(0)$ Levi-flat. Its complexification can be written as $$F_{\mathbb{C}}(x,y,z,w)=\frac{1}{2}P(x,y)+\frac{1}{2}P(z,w)+H_{\mathbb{C}}(x,y,z,w).$$
\par We take the exceptional divisor $D=\mathbb{P}^{3}$ of the blow-up $\pi:(\tilde{\mathbb{C}}^{4},\mathbb{P}^{3})\rightarrow(\mathbb{C}^{4},0)$ with homogeneous coordinates $[a:b:c:d]$, $(a,b,c,d)\in\mathbb{C}^{4}\backslash\{0\}$. The intersection of the strict transform $\tilde{M}_{\mathbb{C}}$ of $M_{\mathbb{C}}$ by $\pi$ with the divisor $D=\mathbb{P}^{3}$ is the surface $$Q=\{[a:b:c:d]\in\mathbb{P}^{3}/P(a,b)+P(c,d)=0\},$$
which is an irreducible smooth surface.
\par Consider for instance the chart $(W,(t,u,z,v))$ of $\tilde{\mathbb{C}}^{4}$ where $$\pi(t,u,z,v)=(t.z,u.z,z,v.z)=(x,y,z,w).$$ We have 
$$F_{\mathbb{C}}\circ\pi(t,u,z,v)=z^{k}(\frac{1}{2}P(t,u)+\frac{1}{2}P(1,v)+zH_{1}(t,u,z,v)),$$
where $H_{1}(t,u,z,v)=H(tz,uz,z,vz)/z^{k+1}$, which implies that $$\tilde{M}_{\mathbb{C}}\cap W=(\frac{1}{2}P(t,u)+\frac{1}{2}P(1,v)+zH_{1}(t,u,z,v)=0)$$
and so $Q\cap W=(z=P(t,u)+P(1,v)=0)$. \par
On the other hand, as we have seen in $(3.2)$, the foliation $\mathcal{L}_{\mathbb{C}}$ is defined by $\alpha|_{M^{*}_{\mathbb{C}}}=0$, where 
$$\alpha=\frac{1}{2}\frac{\partial{P}}{\partial{x}}dx+\frac{1}{2}\frac{\partial{P}}{\partial{y}}dy+
\frac{\partial{H}_{\mathbb{C}}}{\partial{x}}dx+\frac{\partial{H}_{\mathbb{C}}}{\partial{y}}dy.$$ 

In particular, we get
$$\pi^{*}(\alpha)=z^{k-1}(\frac{1}{2}\frac{\partial{P}}{\partial{x}}(t,u)zdt+\frac{1}{2}\frac{\partial{P}}{\partial{y}}(t,u)zdu+
 \frac{1}{2}kP(t,u)dz+z\theta),$$
where $\theta=\pi^{*}(\frac{\partial{H}_{\mathbb{C}}}{\partial{x}}dx+\frac{\partial{H}_{\mathbb{C}}}{\partial{y}}dy)/z^{k}$.
\par Hence, $\tilde{\mathcal{L}}_{\mathbb{C}}$ is defined by 
\begin{eqnarray}\label{equa-integral}
\alpha_{1}=\frac{1}{2}\frac{\partial{P}}{\partial{x}}(t,u)zdt+\frac{1}{2}\frac{\partial{P}}{\partial{y}}(t,u)zdu+
 \frac{1}{2}kP(t,u)dz+z\theta.
\end{eqnarray}
\par Since $Q\cap W=(z=P(t,u)+P(1,v)=0)$, we see that $Q$ is $\tilde{\mathcal{L}}_{\mathbb{C}}$-invariant. In particular, $S:=Q\backslash sing(\tilde{\mathcal{L}}_{\mathbb{C}})$ is a leaf of $\tilde{\mathcal{L}}_{\mathbb{C}}$. Fix $p_{0}\in S$ and a transverse section $\sum$ through $p_{0}$. Let $G\subset Diff(\sum,p_{0})$ be the holonomy group of the leaf $S$ of $\tilde{\mathcal{L}}_{\mathbb{C}}$.
Since $dim(\sum)=1$, we can think that $G\subset Diff(\mathbb{C},0)$. Let us prove that $G$ is finite and linearizable.
\par At this part we use that the leaves of $\tilde{\mathcal{L}}_{\mathbb{C}}$ are closed (see lemma \ref{closed-lema}). 
\par Let $G'=\{f'(0)/f\in G\}$ and consider the homomorphism $\phi:G\rightarrow G'$ defined by $\phi(f)=f'(0)$. We assert that $\phi$ is injective. In fact, assume that $\phi(f)=1$ and by contradiction that $f\neq id$. In this case $f(z)=z+a.z^{r+1}+\ldots$, where $a\neq 0$. According to [L], the pseudo-orbits of
this transformation accumulate at $0\in(\sum,0)$, contradicting that the leaves of $\tilde{\mathcal{L}}_{\mathbb{C}}$ are closed. Now, it suffices to prove that any element $g\in G$ has finite order (cf. [M-M]). In fact, if $\phi(g)=g'(0)$ is a root of unity then $g$ has finite order because $\phi$ is injective. On the other hand, if $g'(0)$ was not a root of unity then $g$ would have pseudo-orbits accumulating at $0\in(\sum,0)$ (cf. [L]). Hence, all transformations of $G$ have finite order and $G$ is linearizable. 
\par This implies that there is a coordinate system $w$ on $(\sum,0)$ such that $G=\langle w\rightarrow\lambda w\rangle$, where $\lambda$ is a $d^{th}$-primitive root of unity (cf. [M-M]). In particular, $\psi(w)=w^{d}$ is a first integral of $G$, that is $\psi\circ g=\psi$ for any $g\in G$.
\par Let $Z$ be the union of the separatrices of $\mathcal{L}_{\mathbb{C}}$ through $0\in\mathbb{C}^{4}$ and $\tilde{Z}$ be its strict transform under $\pi$. The first integral $\psi$ can be extended to a first integral $\varphi:\tilde{M}_{\mathbb{C}}\backslash\tilde{Z}\rightarrow\mathbb{C}$ be setting $$\varphi(p)=\psi(\tilde{L}_{p}\cap \sum),$$ where $\tilde{L}_{p}$ denotes the leaf of $\tilde{\mathcal{L}}_{\mathbb{C}}$ through $p$. Since $\psi$ is bounded (in a compact neighborhood of $0\in\sum$), so is $\varphi$. It follows from Riemann extension theorem that $\varphi$ can be extended holomorphically to $\tilde{Z}$ with $\varphi(\tilde{Z})=0$. This provides the first integral and finishes the proof of theorem \ref{main-theorem}.

\section{Quasihomogeneous polynomials}
In this section, we state some general facts about normal forms of quasihomogeneous polynomials. 
\begin{definition}  
 The Newton support of germ $f=\sum a_{ij}x^{i}y^{j}$ is defined as 
$supp(f)=\{(i,j): a_{ij}\neq 0\}$.
\end{definition}

\begin{definition}
 A holomorphic function $f:(\mathbb{C}^{n},0)\rightarrow(\mathbb{C},0)$ is said to be quasihomogeneous of degree $d$ with indices $\alpha_{1},\ldots,\alpha_{n}$, if for any $\lambda\in\mathbb{C}$ and $(z_{1},\ldots,z_{n})\in\mathbb{C}^{n}$, we have $$f(\lambda^{\alpha_{1}}z_{1},\ldots,\lambda^{\alpha_{n}}z_{n})=\lambda^{d}f(z_{1},\ldots,z_{n}).$$ The index $\alpha_{s}$ is also called the weight of the variable $z_{s}$.
\end{definition}
In the above situation, if $f=\sum a_{k}x^{k}$, $k=(k_{1},\ldots,k_{n})$, $x^{k}=x_{1}^{k_{1}}\ldots x^{k_{n}}$, then $supp(f)\subset\Gamma=\{k:a_{1}k_{1}+\ldots+a_{n}k_{n}=d\}$. The set $\Gamma$ is called the diagonal. Usually one takes $\alpha_{i}\in\mathbb{Q}$ and $d=1$.

One can define the quasihomogeneous filtration of the ring $\mathcal{O}_{n}$. It consists of the decreasing family of ideals $\mathcal{A}_{d}\subset\mathcal{O}_{n}$, $\mathcal{A}_{d'}\subset\mathcal{A}_{d}$ for $d<d'$. Here $\mathcal{A}_{d}=\{Q:$ degrees of monomials from $supp(Q)$ are $deg(Q)\geq d\}$; (the degree is quasihomogeneous). 

When $\alpha_{1}=\ldots=\alpha_{n}=1$, this filtration coincides with the usual filtration by the usual degree.

\begin{definition}
A function $f$ is called semiquasihomogeneous if $f=Q+F'$, where $Q$ is quasihomogeneous of degree $d$ of finite multiplicity and $F'\in\mathcal{A}_{d'}$, $d'>d$. 
\end{definition}
\par We will use the following result (cf. [A]). 
\begin{theorem}~\label{arno-teo}
Let $f$ be a semiquasihomogeneous function, $f=Q+F'$ with quasihomogeneous $Q$ of finite multiplicity. Then $f$ is right equivalent to the function $Q+\sum_{j} c_{j}e_{j}(z)$, where $e_{1},\ldots,e_{s}$ are the elements of the monomial basis of the local algebra $A_{Q}$ such that $deg(e_{j})>d$ and $c_{j}\in\mathbb{C}$. 
\end{theorem}
\begin{example}
 If $f=Q+F'$ is semiquasihomogeneous and $Q(x,y)=x^{2}y+y^{k}$, then $f$ is right equivalent to $Q$. Indeed, the base of the local algebra $\mathcal{O}_{2}/(xy,x^{2}+ky^{k-1})$  is $1,x,y,y^{2},\ldots,y^{k-1}$ and lies below the diagonal $\Gamma$. Here $\mu(Q,0)=k+1$.
\end{example}

\section{Proof of theorem \ref{second-theorem}}
Let $M=F^{-1}(0)$ be a germ at $0\in\mathbb{C}^{n}$, $n\geq 3$ of real analytic Levi-flat hypersurface, where $F(z)=\mathcal{R}e(Q(z))+h.o.t$ and $Q$ is a quasihomogeneous polynomial of degree $d$ with an isolated singularity at $0\in\mathbb{C}^{n}$. It is easily seen that $sing(M_{\mathbb{C}})=\{0\}$ and $cod_{M^{*}_{\mathbb{C}}}sing(\mathcal{L}_{\mathbb{C}})\geq 3$. The argument is essentially the same of the proof of theorem \ref{main-theorem}. In this way, there exists an unique germ at $0\in\mathbb{C}^{n}$ of holomorphic codimension one foliation $\mathcal{F}_{M}$ tangent to $M$, moreover $\mathcal{F}_{M}$: $dh=0$, $h(z)=Q(z)+h.o.t$ and $M=(\mathcal{R}e(h)=0)$. Acoording to theorem \ref{arno-teo}, there exists $\phi\in Diff(\mathbb{C}^{n},0)$ such that $h\circ \phi^{-1}(w)=Q(w)+\sum_{k}c_{k}e_{k}(w)$, where $c_{k}$ and $e_{k}$ as above. Hence $$\phi(M)=(\mathcal{R}e(Q(w)+\sum_{k}c_{k}e_{k}(w))=0).$$

\section{Applications}
\par Here we give some applications of theorem \ref{main-theorem}.
\begin{example}\label{exa-arturo}
 $Q(x,y)=x^{2}y+y^{3}$ is a homogeneous polynomial of degree $3$ with an isolated singularity at $0\in\mathbb{C}^{2}$ and Milnor number $\mu(Q,0)=4$. According to [A-G-V] pg. 184, any germ $f(x,y)=x^{2}y+y^{3}+h.o.t$ is right equivalent to $x^{2}y+y^{3}$.
\par In particular, if $F(z)=\mathcal{R}e(x^{2}y+y^{3}) + h.o.t$ and $M=(F=0)$ is a germ of real analytic Levi-flat at $0\in\mathbb{C}^{2}$, theorem \ref{main-theorem} implies that there exists a holomorphic change of coordinate such that $$M=(\mathcal{R}e(x^{2}y+y^{3})=0).$$
\end{example}
\begin{example}

\par If $Q(x,y)=x^{5}+y^{5}$ then $f(x,y)=Q(x,y)+h.o.t$ is right equivalent to $x^{5}+y^{5}+c.x^{3}y^{3}$, where $c\neq 0$ is a constant (cf. [A-G-V], pg. 194). Let $F(z)=\mathcal{R}e(x^{5}+y^{5}) + h.o.t$ be such that $M=(F=0)$ is Levi-flat, theorem \ref{main-theorem} implies that there exists a holomorphic change of coordinate such that $$M=(\mathcal{R}e(x^{5}+y^{5}+c.x^{3}y^{3})=0).$$
 
\end{example}

\begin{center}
 \noindent{\large\bf References}
\end{center}

\begin{itemize}
\item[{[A]}] V.I. Arnold:
 ``Normal Form of functions in the neighbourhood of degenerate critical points", 
 UNM 29:2 (1974), 11-49, RMS 29:2 19-48.
\item[{[A-G-V]}]  V.I. Arnold, S.M. Gusein-Zade, A.N. Varchenko:
   ``Singularities of Differential Maps",
    Vol. I, Monographs in Math., vol. 82, Birkh\"{a}user, $1985$.
\item[{[B-G]}]  D. Burns, X. Gong:
   ``Singular Levi-flat real analytic hypersurfaces",
    Amer. J. Math. 121, $(1999)$, pp. 23-53.

\item[{[C-L]}]  D. Cerveau, A. Lins Neto:
   ``Local Levi-Flat hypersurfaces invariants by
             a codimension one holomorphic foliation".
     To appear in Amer. J. Math.

\item[{[L]}] F. Loray: ``Pseudo-groupe d'une singularit\' e de feuilletage holomorphe en dimension deux". Avaliable in \\
http://hal.archives-ouvertures.fr/ccsd-00016434

\item[{[M-M]}] J.F. Mattei, R. Moussu:
   ``Holonomie et int\' egrales premi\` eres",
    Ann. Ec. Norm. Sup. 13, $(1980)$, pg. 469-523.

\end{itemize} 
 \end{document}